\documentclass[12pt]{amsart}

\usepackage{amsmath}
\usepackage{amssymb}

\def\Gal{\mathrm{Gal}}
\def\zprank{\mathrm{rank}_{\mathbb{Z}_p}}
\def\charLambda{\mathrm{char}_{\Lambda}}

\title{On tamely ramified Iwasawa modules for the cyclotomic $\mathbb{Z}_p$-extension 
of abelian fields
\footnote{2010 Mathematics Subject Classification. Primary 11R23, 
Secondary 11R18}}
\author{Tsuyoshi Itoh}
\date{December 6, 2012.}

\begin{document}

\maketitle

\begin{abstract}
Let $p$ be an odd prime, and 
$k_\infty$ the cyclotomic $\mathbb{Z}_p$-extension of an abelian field $k$.
For a finite set $S$ of rational primes which does not include $p$,
we will consider the maximal $S$-ramified abelian pro-$p$ extension 
$M_S (k_\infty)$ over $k_\infty$.
We shall give a formula of the $\mathbb{Z}_p$-rank of 
$\Gal (M_S (k_\infty)/k_\infty)$.
In the proof of this formula, we also show that 
$M_{\{q\}} (k_\infty)/L(k_\infty)$ is a finite extension
for every real abelian field $k$ and every rational prime $q$ distinct from $p$,
where $L(k_\infty)$ is the maximal unramified abelian pro-$p$ extension over $k_\infty$.
\end{abstract}

\section{Introduction}
Let $k$ be an algebraic number field, and $p$ a prime number.
We denote by $k_\infty/k$ the cyclotomic $\mathbb{Z}_p$-extension 
(i.e. the unique $\mathbb{Z}_p$-extension contained in the field generated by
all $p$-power roots of unity over $k$).
Let $S$ be a finite set of \textit{rational} primes,
and $\widetilde{M}_S (k_\infty)$ the maximal pro-$p$ extension of $k_\infty$ 
unramified outside $S$
(i.e., the primes of $k_\infty$ lying above the primes in $S$
are only allowed to ramify in $\widetilde{M}_S (k_\infty)/k_\infty$).

When $p \in S$, the structure of 
$\widetilde{X}_S (k_\infty) = \Gal (\widetilde{M}_S (k_\infty)/k_\infty)$ 
is already studied (see, e.g., Iwasawa \cite{Iwa81}, Neukirch-Schmidt-Wingberg \cite{NSW}).
In particular, $\widetilde{X}_S (k_\infty)$ is a free pro-$p$ group under 
certain conditions.

Recently, 
the structure of $\widetilde{X}_S (k_\infty)$ 
for the case that $p \not\in S$ is also studied
by several authors (Salle \cite{Salle}, Mizusawa-Ozaki \cite{MO}, ...).
In this case, it seems that $\widetilde{X}_S (k_\infty)$ does not have a simple structure.
Then, to study $\widetilde{X}_S (k_\infty)$, it is important to
study the structure of its abelian quotient.
Let $M_S (k_\infty)/k_\infty$ be the maximal abelian pro-$p$ extension
unramified outside $S$.
In the present paper, we shall consider 
$X_S (k_\infty) = \Gal (M_S (k_\infty)/k_\infty)$ 
for the case that $p \not\in S$.
Since $\Gal (k_\infty/k)$ acts on $X_S (k_\infty)$, we can use Iwasawa
theoretic arguments.
We call this $X_S (k_\infty)$ the $S$-ramified Iwasawa module.
(When $p \not\in S$, all primes which ramify in $M_S (k_\infty)/k_\infty$
are tamely ramified.
We also call these modules ``tamely ramified Iwasawa modules''.)

If a $\mathbb{Z}_p$-module $M$ satisfies 
$\dim_{\mathbb{Q}_p} M \otimes_{\mathbb{Z}_p} \mathbb{Q}_p = r < \infty$, 
we say that the $\mathbb{Z}_p$-rank of $M$ is $r$, and we write $\zprank M =r$.
Our purpose of the present paper is giving a formula of
$\zprank X_S (k_\infty)$ when $k$ is an abelian extension 
of $\mathbb{Q}$ (abelian field) and $p$ is an odd prime.
(In this case, we can show that $X_S (k_\infty)$ is finitely generated over $\mathbb{Z}_p$.)

We shall give a remark about ``giving a formula of $\zprank X_S (k_\infty)$''.
Let $L(k_\infty)$ be the maximal unramified abelian pro-$p$ extension of $k_\infty$.
We put $X(k_\infty)=\Gal (L(k_\infty)/k_\infty)$.
Then $X(k_\infty)$ is the (usual) Iwasawa module,
and $\lambda = \zprank X(k_\infty)$ is called the Iwasawa $\lambda$-invariant.
In general, it is hard to write $\lambda$ explicitly.
Since $X(k_\infty)$ is a quotient of $X_S (k_\infty)$,
we consider it is sufficient to obtain a formula including $\lambda$ 
at the present time.
(That is, we will only give a formula of $\zprank \Gal (M_S (k_\infty)/L (k_\infty))$,
actually.)
However, for abelian fields, the ``plus part'' of $\lambda$ is conjectured to be $0$
(Greenberg's conjecture),
and the ``minus part'' of $\lambda$ can be computed (at least theoretically)
from the Kubota-Leopoldt $p$-adic $L$-functions
(Stickelberger elements).

We also mention that 
formulas of $\zprank X_S (k_\infty)$ are already obtained for several cases.
In particular, it can be said that the $\mathbb{Z}_p$-rank of the
``minus part'' of $X_S (k_\infty)$ for 
CM-fields is already known (see section 2).
Salle \cite{Salle} studied $X_S (k_\infty)$
for the case that $k$ is an imaginary quadratic field (or $\mathbb{Q}$) with $p=2$.
Moreover, when $k=\mathbb{Q}$, 
a formula of $\zprank X_S (\mathbb{Q}_\infty)$ (including the case that $p=2$) 
is shown by Mizusawa, Ozaki, and the author \cite{IMO}
(as a corollary, a general formula for imaginary quadratic fields is also given).
In the present paper, we shall extend the method given in \cite{IMO} for 
abelian fields.
The following theorem is crucial
to prove the formula of $\zprank X_S (\mathbb{Q}_\infty)$.

\par \bigskip

\noindent
\textbf{Theorem A} (see \cite{IMO}). 
{\itshape
Let $q$ be a rational prime distinct from $p$. 
Then $M_{\{q\}} (\mathbb{Q}_\infty)/\mathbb{Q}_\infty$ is a finite extension.
}

\par \bigskip

At first, we will generalize Theorem A for real abelian fields
(under the condition that $p$ is an odd prime).

\par \bigskip

\noindent
\textbf{Theorem 1.1. } 
{\itshape
Assume that $p$ is odd. 
Let $k$ be a real abelian field, 
and $q$ a rational prime distinct from $p$. 
Then $M_{\{q\}} (k_\infty) / L(k_\infty)$ is a finite extension.
}

\par \bigskip

We remark that $M_{\{q\}} (k_\infty) / L(k_\infty)$ can be infinite
when $k$ is an imaginary abelian filed.
(For example, see \cite{Salle}, \cite{KW}, \cite{IMO}, 
or section 6 of the present paper.)
Similar to \cite{IMO}, Theorem 1.1 plays an important role to prove 
our formula of $\zprank X_S (k_\infty)$ for abelian fields.

In section 2, 
we shall state some basic facts, and give preparations for proving Theorem 1.1.
We will prove Theorem 1.1 in sections 3 and 4.
In section 5, we shall give a simple remark about a generalization 
of Theorem 1.1.
In section 6, 
we shall give a formula of $\zprank X_S (K_\infty)$ for abelian fields
(Theorem 6.4). 
The formula is given as the $\chi$-component version.
We also give examples with applying this formula
for some simple cases.

\section{Preliminaries}
Firstly, we shall recall some basic facts from class field theory.
Let $p$ be an odd prime number, and $k$ an algebraic number field.
(In the following of the present paper, we assume that $p$ is odd.)
We denote by $k_\infty/k$ the cyclotomic $\mathbb{Z}_p$-extension.
For a non-negative integer $n$, let $k_n$ be the $n$th layer of $k_\infty/k$
(that is, the unique subfield of $k_\infty$ such that $k_n/k$ is a cyclic extension 
of degree $p^n$).
Let $S$ be a finite set of rational primes which does not include $p$.
For an algebraic extension (not necessary finite) $\mathcal{K}$ of $\mathbb{Q}$, 
let $M_S (\mathcal{K})$ be the maximal abelian (pro-)$p$-extension of 
$\mathcal{K}$ unramified outside $S$,
and $L(\mathcal{K})$ the maximal unramified abelian (pro-)$p$-extension of $\mathcal{K}$.
For an abelian group $G$, let $\widehat{G}$ be the $p$-adic completion of $G$
(that is, $\widehat{G}=\varprojlim G/G^{p^n}$).

As noted in section 1, 
we shall mainly consider the $\mathbb{Z}_p$-rank of  
$\Gal (M_S (k_\infty)/L(k_\infty))$.
We will write several facts which is also stated in \cite{IMO}.
In this paragraph, assume that $S$ is not empty.
By class field theory, we have the following exact sequence:
\[ \widehat{E_{k_n}} \overset{\eta_n}{\rightarrow} \bigoplus_{q \in S} 
\widehat{(O_{k_n}/q)^{\times}} \rightarrow \Gal (M_S (k_n)/L(k_n)) \rightarrow 0, \]
where $E_{k_n}$ is the group of units of $k_n$, $O_{k_n}$ is the 
ring of integers of $k_n$, 
and $\eta_n$ is the natural homomorphism induced from the diagonal embedding.
(We will give a remark on the structure of $\widehat{(O_{k_n}/q)^{\times}}$.
Assume that the prime decomposition of $q O_{k_n}$ is 
$\mathfrak{q}_1^{e_1} \cdots \mathfrak{q}_r^{e_r}$.
Then 
\[ \widehat{(O_{k_n}/q)^{\times}} \cong \bigoplus_{i=1}^r 
\widehat{(O_{k_n}/\mathfrak{q}_i)^{\times}} \]
because $\widehat{(O_{k_n}/\mathfrak{q}_i^{e_i})^{\times}} \cong 
\widehat{(O_{k_n}/\mathfrak{q}_i)^{\times}}$.)
We put $E_\infty = \varprojlim \widehat{E_{k_n}}$, and 
$R_q = \varprojlim \widehat{(O_{k_n}/q)^{\times}}$, where the projective limits
are taken with respect to the natural mappings induced from the norm mapping.
Then we obtain the following exact sequence:
\[ E_\infty \overset{\eta_\infty}{\rightarrow} \bigoplus_{q \in S} R_q \rightarrow 
\Gal (M_S (k_\infty)/L(k_\infty)) \rightarrow 0. \]
In the cyclotomic $\mathbb{Z}_p$-extension, all (finite) primes of $k$ are finitely decomposed.
Then $R_q$ is a finitely generated $\mathbb{Z}_p$-module.
From this, we also see $\Gal (M_S (k_\infty)/L(k_\infty))$ is 
finitely generated over $\mathbb{Z}_p$.
On the other hand, the theorem of Ferrero-Washington \cite{FW} implies that
$\Gal (L(k_\infty)/k_\infty)$ is a finitely generated $\mathbb{Z}_p$-module,
if $k$ is an abelian field.
Hence, we see that for every abelian field $k$, 
$X_S (k_\infty)$ is finitely generated over $\mathbb{Z}_p$.
(We can see that the $\mathbb{Z}_p$-rank of $X_S (k_\infty)$ is always finite
in general.)
It seems hard to determine the cokernel of $\eta_\infty$ directly.

\par \bigskip

\noindent
\textbf{Remark. }
When the base field $k$ is a CM-field, the minus part of $E_\infty$ is 
easy to compute (we are also able to compute the minus part of $R_q$).
Hence we can obtain a formula of the $\mathbb{Z}_p$-rank of the minus part 
of $\Gal (M_S (k_\infty)/L(k_\infty))$.
This idea is already known (see, e.g., \cite{Salle}, \cite{KW}, \cite{IMO}).

\par \bigskip

Secondly, we shall give some preparations to prove Theorem 1.1.
Let $q$ be a prime number satisfying $q \neq p$.
For simplicity, we will write $M_p ( \cdot )$, $M_q ( \cdot )$ instead of
$M_{\{p\}} ( \cdot )$, $M_{\{q\}} ( \cdot )$, respectively.

\par \bigskip

\noindent
\textbf{Lemma 2.1. } 
{\itshape
Let $k'/k$ be a finite extension of algebraic number fields.
If $\Gal(M_q (k'_\infty) / L(k'_\infty))$ is finite, 
then $\Gal(M_q (k_\infty) / L(k_\infty))$ is also finite.
}

\par \bigskip

\noindent
\textbf{Proof. } 
We may assume that $k' \cap k_\infty=k$.
Let 
\[ N_n : \widehat{(O_{k'_n}/q)^{\times}} \rightarrow \widehat{(O_{k_n}/q)^{\times}} \]
be the homomorphism induced from the norm mapping.
We can see that the order of the cokernel $\mathrm{Coker} (N_n)$ is bounded as $n \to \infty$.
(Proof: Since there are only finitely many primes in $k_\infty$ lying above $q$,
the $p$-rank of $\widehat{(O_{k_n}/q)^{\times}}$ is bounded.
Moreover, the exponent of $\mathrm{Coker} (N_n)$ is at most $[k':k]$.)
Since $\Gal(M_q (k'_n) / L(k'_n))$ (resp. $\Gal(M_q (k_n) / L(k_n))$)
is isomorphic to a quotient of 
$\widehat{(O_{k'_n}/q)^{\times}}$ (resp. $\widehat{(O_{k_n}/q)^{\times}}$),
$N_n$ induces the homomorphism 
$\Gal(M_q (k'_n) / L(k'_n)) \rightarrow \Gal(M_q (k_n) / L(k_n))$.
From the above fact, the order of the cokernel is bounded as $n \to \infty$.

Assume that $\Gal(M_q (k'_\infty) / L(k'_\infty))$ is finite.
Then we can show that the order of $\Gal(M_q (k'_n) / L(k'_n))$ is bounded as $n \to \infty$.
From the above fact, we see that the order of $\Gal(M_q (k'_n) / L(k'_n))$ 
is also bounded.
Hence $\Gal(M_q (k_\infty) / L(k_\infty))$ is finite.
\hfill $\Box$

\par \bigskip

From Lemma 2.1 and the theorem of Kronecker-Weber, 
we may replace a real abelian field $k$ 
to the maximal real subfield of a cyclotomic filed containing $k$ 
to show Theorem 1.1.
For a positive integer $d$,
let $\mu_d$ be the set of all $d$th roots of unity, and $\mathbb{Q} (\mu_d)$ the 
$d$th cyclotomic field.

\par \bigskip

\noindent
\textbf{Lemma 2.2. } 
{\itshape
Let $f$ be a positive integer which is prime to $p$, and $m$ a positive integer.
We put $K=\mathbb{Q} (\mu_{fp^m})$ and $k=K^+$ 
(the maximal real subfield of $K$).
If $q$ does not split in $K/k$, then $M_q (k_\infty)=L(k_\infty)$.
}

\par \bigskip

\noindent
\textbf{Proof. } 
Let $\mathfrak{q}$ be an arbitrary prime of $k$ lying above $q$.
It is well known that if $\mathfrak{q}$ does not split in $K$, then 
the order of $(O_k/\mathfrak{q})^{\times}$ is not divisible by $p$.
(Proof: We denote by $k_\mathfrak{q}$ the completion of $k$ at $\mathfrak{q}$.
Under the assumption, $k_\mathfrak{q}$ does not contain $\mu_p$.
By the structure of the group of units in $k_\mathfrak{q}$, 
we obtain the assertion.)
Since $\Gal(M_q (k)/L(k))$ is isomorphic to a quotient of $(O_k/\mathfrak{q})^{\times}$,
we see $M_q (k)=L(k)$.

We note that the $n$th layer $k_n$ of $k_\infty/k$ is the maximal real 
subfield of $\mathbb{Q} (\mu_{fp^{m+n}})$.
Hence by using the same argument, we also see $M_q (k_n)=L(k_n)$
for all $n \geq 1$.
This implies that $M_q (k_\infty)=L(k_\infty)$.
\hfill $\Box$

\par \bigskip

From the above arguments, it is sufficient to prove Theorem 1.1 under the following 
conditions:

\par \bigskip

\noindent \textbf{(A)}
$k$ is the maximal real subfield of $K=\mathbb{Q} (\mu_{fp^m})$,
where $f$ and $m$ are positive integers and $f$ is prime to $p$. 
Every prime lying above $q$ splits in $K/k$,
and is not decomposed in $k_\infty/k$
(the latter can be satisfied by taking $m$ sufficiently large).

\section{Properties of certain Kummer extensions}
In this section, we shall give some key results to prove Theorem 1.1.
Assume that $K, k$, and $q$ satisfy (A) in section 2.
We will construct certain infinite Kummer extensions over $K_\infty$.
We shall use some fundamental results given by Khare-Wintenberger \cite{KW}.

We define the terms \textbf{case NS} and \textbf{case S} as follows:

\par \bigskip

\noindent \begin{tabular}{lcl}
case NS & : & every prime lying above $p$ does not split in $K/k$, \\
case S  & : & every prime lying above $p$ splits in $K/k$.
\end{tabular}

\par \bigskip

Moreover, we use the following notation (in sections 3 and 4):

\begin{itemize}
\item $J$ : complex conjugation

\item $\mathfrak{q}_1, \ldots, \mathfrak{q}_r$ : prime ideals of $k$ lying above $q$

\item $\mathfrak{p}_1, \ldots, \mathfrak{p}_t$ : prime ideals of $k$ lying above $p$

\item $\mathfrak{Q}_i, \mathfrak{Q}^J_i$ $(i=1, \ldots, r)$ :
prime ideals of $K$ lying above $\mathfrak{q}_i$

\item $\mathfrak{P}_j$ $(j=1, \ldots, t)$ : (unique) prime ideal of $K$ lying above 
$\mathfrak{p}_j$ (case NS)

\item $\mathfrak{P}_j, \mathfrak{P}^J_j$ $(j=1, \ldots, t)$ : 
prime ideals of $K$ lying above $\mathfrak{p}_j$ (case S)
\end{itemize}

Following Greenberg \cite{Gre73}, we denote by $s$ the number of primes of $k$
which is lying above $p$ and splits in $K$.
Hence we see that $s=0$ for the case NS, and $s=t$ for the case S.
Note that every prime lying above $p$ are totally ramified in $k_\infty/k$
by the assumption on $k$.
Hence $s$ is also the number of primes of $k_\infty$
which is lying above $p$ and splits in $K_\infty$.

By the assumption, $K_\infty$ contains all $p^n$th roots of unity.
For an element $x$ of $K^{\times}$, we define
\[ K_\infty (\sqrt[p^\infty]{x}) 
= \bigcup_{n \geq 1} K_\infty (\sqrt[p^n]{x}). \]
More precisely, $K_\infty (\sqrt[p^\infty]{x})$ is 
the union of all finite Kummer extensions 
$K_n (\sqrt[p^n]{x})$ for $n \geq 1$
(note that $K_n$ contains $\mu_{p^n}$).
Similarly, for a finitely generated subgroup $T$ of $K^{\times}$, 
we define the extension 
$K_\infty (\sqrt[p^\infty]{T})/ K_\infty$
by adjoining all $p^n$th roots of the elements contained in $T$.
As noted in \cite{KW}, $K_\infty (\sqrt[p^\infty]{x}) =K_\infty$ if and only if 
$x$ is a root of unity.

The following result is helpful to prove the results stated in this section.

\par \bigskip

\noindent
\textbf{Theorem B} (see Khare-Wintenberger \cite[Lemma 2.5]{KW}). 
{\itshape
Let $T$ be a finitely generated subgroup of $K^{\times}$, 
and $S$ a finite set of (finite) primes of $K$.
Let $\mathcal{I}$ be the subgroup of $\Gal(K_\infty (\sqrt[p^\infty]{T})/ K_\infty)$
generated by the inertia subgroups for the primes in $S$.
For a prime $\mathfrak{r} \in S$, let $K_\mathfrak{r}$ be the 
completion of $K$ at $\mathfrak{r}$.
We denote by $\mathcal{T}$ the closure of the diagonal image of $T$ in 
$\prod_{\mathfrak{r} \in S} \widehat{K_\mathfrak{r}^{\times}}$.
(Recall that $\widehat{K_\mathfrak{r}^{\times}}$ is the $p$-adic completion of 
$K_\mathfrak{r}^{\times}$.)
Then $\zprank \mathcal{I} = \zprank \mathcal{T}$. 
}

\par \bigskip

We shall construct several Kummer extensions unramified outside $\{p,q\}$ over $K_\infty$
by following the method given in \cite{IMO}.
(See also Greenberg \cite{Gre77}.)
Let $k^D$ be the decomposition field of $K/\mathbb{Q}$ for $q$.
By the assumption, $k^D$ is an imaginary abelian field 
and $[ k^D : \mathbb{Q} ] =2r$.
Let $Q_1, \ldots, Q_r, Q_1^J, \ldots, Q_r^J$ be the primes of $k^D$
lying below $\mathfrak{Q}_1, \ldots, \mathfrak{Q}_r, 
\mathfrak{Q}_1^J, \ldots, \mathfrak{Q}_r^J$ respectively.
We can take a positive integer $h$ such that 
\begin{itemize}
\item $Q_1^h$ is a principal ideal generated by $\alpha_1$, and 
\item $\alpha_1-1 \in P$ for every prime ideal $P$ of $k^D$ lying above $p$.
\end{itemize}
We note that $Q_1^\sigma$ $(\sigma \in \Gal(k^D/\mathbb{Q}))$ is the complete set of 
primes in $k^D$ lying above $q$, and $(Q_1^\sigma)^h = (\alpha_1^\sigma)$.
We write all conjugates of $\alpha_1$ for $\Gal(k^D/\mathbb{Q})$ as the following:
\[ \alpha_1, \alpha_2, \ldots, \alpha_r, 
\alpha_1^J, \alpha_2^J, \ldots, \alpha_r^J \]
(these are distinct elements because $q$ splits completely in $k^D$).
Moreover, we put $\beta_i =\alpha_i/\alpha_i^J$ for $i=1, \ldots, r$.

For the case S (i.e. $p$ splits in $K/k$), 
we define $\rho_1, \ldots, \rho_t \in K^{\times}$ as follows.
We can take an integer $h'$ such that $\mathfrak{P}_j^{h'} = (\pi_j)$ 
for all $j=1, \ldots, t$.
We put $\rho_j =\pi_j/\pi_j^J$ for $j=1, \ldots, t$.

\par \bigskip

\noindent
\textbf{Definition. } 
We put 
\[ T_q = \langle \beta_i \; | \; i=1, \ldots, r \rangle, \]
which is a subgroup of $(k^D)^{\times}$ (and hence also a 
subgroup of $K^{\times}$).
For the case NS, we put $T=T_q$.
For the case S, we put
\[ T = \langle \beta_i, \rho_j \; | \; i=1, \ldots, r, \; j=1, \ldots, t \rangle. \]
Moreover, we put $N_q=K_\infty (\sqrt[p^\infty]{T_q})$ and 
$N=K_\infty (\sqrt[p^\infty]{T})$.
(Of course, $N=N_q$ for the case NS.)

\par \bigskip

\noindent
\textbf{Lemma 3.1. }
(1) {\itshape $N$ and $N_q$ are abelian extensions over $k_\infty$.}
(2) {\itshape $N/K_\infty$ and $N_q/K_\infty$ are unramified outside $\{p,q\}$.}
(3) {\itshape $\Gal(N/K_\infty) \cong \mathbb{Z}_p^{\oplus r+s}$ and
$\Gal(N_q/K_\infty) \cong \mathbb{Z}_p^{\oplus r}$.
(Recall that $s=0$ for the case NS, and $s=t$ for the case S.)
}

\par \bigskip

\noindent
\textbf{Proof. }
(1) Since $J$ acts on $T$ as $-1$,
then $J$ acts on $\Gal(N/K_\infty)$ and the action is trivial.
This implies that $N/k_\infty$ is an abelian extension.
The assertion for $N_q$ follows similarly.

(2) Note that all elements contained in $T$ (resp. $T_q$) are $\{p,q\}$-units.
Hence $N/K_\infty$ (resp. $N_q/K_\infty$) is unramified outside $\{p,q\}$.

(3) It is easy to see that $T$ is a free $\mathbb{Z}$-module 
of rank $r+s$. 
Let $\widehat{T}$ be the closure of $T$ in $\widehat{K^{\times}}$.
Then the $\mathbb{Z}_p$-rank of $\widehat{T}$ is $r+s$.
As noted in \cite{KW} (see \textit{Remarks} after the proof of \cite[Lemma 2.2]{KW}),
this fact implies that $\zprank \Gal(N/K_\infty)=r+s$.
Since $\Gal (N/K_\infty)$ is generated by $r+s$ elements, we obtain 
the isomorphism $\Gal(N/K_\infty) \cong \mathbb{Z}_p^{\oplus r+s}$.
The assertion for $\Gal(N_q/K_\infty)$ can be proven quite similarly.
\hfill $\Box$

\par \bigskip

\noindent
\textbf{Lemma 3.2} (cf. Greenberg \cite{Gre77}). \, 
(1) {\itshape For every $i=1, \ldots, r$, the unique prime lying above 
$\mathfrak{Q}_i$ is ramified in $K_\infty (\sqrt[p^\infty]{\beta_i}) / K_\infty$.}
(2) {\itshape $N_q \cap M_p (K_\infty)$ is a finite extension 
over $K_\infty$.}

\par \bigskip

\noindent
\textbf{Proof. }
The assertions can be shown easily by using Theorem B.
Note that (2) is already mentioned in \cite[p.149]{Gre77}.
\hfill $\Box$

\par \bigskip

\noindent
\textbf{Proposition 3.3. }
{\itshape
Let $\mathcal{I}_p$ be the subgroup of $\Gal(N/K_\infty)$ generated by all
inertia groups for the prime lying above $p$.
Then $\mathcal{I}_p$ has finite index in $\Gal(N/K_\infty)$.
}

\par \bigskip

\noindent
\textbf{Proof. }
Firstly, we consider the case NS.
We assume that $p$ does not split in $K/k$.
Hence $s=0$, $T=T_q$, $N=N_q$, and there are just $t$ primes in $K$ lying above $p$.
By Lemma 3.1, 
it is sufficient to show that $\zprank \mathcal{I}_p=r$.
Let $\mathcal{T}_q$ be the closure of the diagonal image of $T_q$ in 
$\prod_{j=1}^t \widehat{K_{\mathfrak{P}_j}^{\times}}$.
By Theorem B, we see that $\zprank \mathcal{I}_p = \zprank \mathcal{T}_q$,
hence we shall show $\zprank \mathcal{T}_p=r$.

Recall that $k^D$ is the decomposition field of $K/\mathbb{Q}$ for $q$,
and $T_q$ is also a subgroup of $(k^D)^{\times}$.
We denote by $P_1, \ldots, P_u$ the primes of $k^D$ lying above $p$
(where $u \leq t$).
Let $\mathcal{T}'_q$ be the closure of the diagonal image of $T_q$ in 
$\prod_{h=1}^u \widehat{(k^D_{P_h})^{\times}}$.

We claim that $\zprank \mathcal{T}'_q =\zprank \mathcal{T}_q$.
By the definition of $T_q$, 
every element $x$ of $T_q$ satisfies $x-1 \in P_h$ for all $h=1, \ldots, u$.
Let $\mathcal{U}^1_{P_h}$ be the group of principal units 
of $k^D_{P_h}$.
We see that $\mathcal{T}'_q$ is contained in $\prod_{h=1}^u \mathcal{U}^1_{P_h}$.
Let $\iota$ be the homomorphism 
\[ \prod_{h=1}^u \mathcal{U}^1_{P_h} \rightarrow 
\prod_{j=1}^t \mathcal{U}^1_{\mathfrak{P}_j} \]
induced from the diagonal embedding 
$\mathcal{U}^1_{P_h} \rightarrow 
\prod_{\mathfrak{P} | p_h} \mathcal{U}^1_{\mathfrak{P}}$.
We can see that $\iota$ is injective, and $\iota (\mathcal{T}'_q)= \mathcal{T}_q$.
Then the claim follows.

We shall recall the argument given in Brumer's proof of Leopoldt's conjecture 
for abelian fields (see also \cite{Brumer}, \cite{Was}).
Assume that $\zprank \mathcal{T}'_q < r$.
Then there are elements $a_1, \ldots, a_r$ of $\mathbb{Z}_p$ which satisfies
\[ \beta_1^{a_1} \beta_2^{a_2} \cdots \beta_r^{a_r} = 1
\quad \mbox{in } \mathcal{U}^1_{P_h} \]
for all $h=1, \ldots, u$, and $a_i \neq 0$ with some $i$.
Since $\beta_i = \alpha_i/\alpha_i^J$ and $\alpha_i$ is a conjugate of $\alpha_1$,
we also see that 
\[ \prod_{\sigma \in \Gal(k^D/\mathbb{Q})} (\alpha_1^{\sigma \tau^{-1}})^{x(\sigma)} = 1
\quad \mbox{in } \mathcal{U}^1_{P_1} \]
for all $\tau \in \Gal(k^D/\mathbb{Q})$, 
where $x(\sigma) \in \mathbb{Z}_p$ 
satisfying $x(\sigma) \neq 0$ with some $\sigma$.
Fix an embedding $k^D_{P_1} \rightarrow \mathbb{C}_p$. 
By taking the (normalized)
$p$-adic logarithm of the above equation, we see
\[ \sum_{\sigma \in \Gal(k^D/\mathbb{Q})} 
x(\sigma) \log_p \alpha_1^{\sigma \tau^{-1}} = 0. \]
This implies that the determinant of the matrix 
$(\log_p \alpha_1^{\sigma \tau^{-1}})_{\sigma, \tau}$ is $0$.

On the other hand, $\alpha_1, \ldots, \alpha_r,
\alpha_1^J, \ldots, \alpha_r^J$ are multiplicative independent in $(k^D)^{\times}$.
Then we can see that $\log_p \alpha_1, \ldots, \log_p \alpha_r,
\log_p \alpha_1^J, \ldots, \log_p \alpha_r^J$ are linearly independent over $\mathbb{Q}$.
By Baker-Brumer's theorem (see Brumer \cite{Brumer}, 
Washington \cite[Theorem 5.29]{Was}),
they are also linearly independent over $\overline{\mathbb{Q}}$ in $\mathbb{C}_p$.
Hence the determinant of the matrix 
$(\log_p \alpha_1^{\sigma \tau^{-1}})_{\sigma, \tau}$ is \textit{not} $0$.
(This follows from the argument given in the proof of \cite[Lemma 6]{IMO}
which uses \cite[Lemma 5.26(a)]{Was}.)
It is a contradiction.
Then we conclude that 
\[ \zprank \mathcal{T}'_q =\zprank \mathcal{T}_q 
= \zprank \mathcal{I}_p = r. \]

Next, we shall consider the case S.
Assume that $p$ splits in $K/k$. 
That is, $s=t$ and the number of primes of $K$ lying above $p$ is $2t$.
The outline of the proof is the same as the case NS.
Let $\mathcal{T}$ be the closure of the diagonal image of $T$ in 
$\prod_{j=1}^t \widehat{K_{\mathfrak{P}_j}^{\times}} \times 
\prod_{j=1}^t \widehat{K_{\mathfrak{P}^J_j}^{\times}}$.
In this case, we shall show $\zprank \mathcal{T} =r+t$.

Assume that $\zprank \mathcal{T} < r+t$.
Then there are elements $a_1, \ldots, a_{r+t}$ of $\mathbb{Z}_p$ which satisfies
\[ \begin{array}{lll}
\beta_1^{a_1} \beta_2^{a_2} \cdots \beta_r^{a_r} 
\rho_1^{a_{r+1}} \rho_2^{a_{r+2}} \cdots \rho_t^{a_{r+t}} = 1
& \mbox{in } \mathcal{U}^1_{\mathfrak{P}_j}, & \mbox{and} \\
\beta_1^{a_1} \beta_2^{a_2} \cdots \beta_r^{a_r} 
\rho_1^{a_{r+1}} \rho_2^{a_{r+2}} \cdots \rho_t^{a_{r+t}} = 1
& \mbox{in } \mathcal{U}^1_{\mathfrak{P}^J_j} & 
\end{array} \]
for all $1 \leq j \leq t$, and $a_i \neq 0$ with some $i$.
However, $v_{\mathfrak{P}_j} (\rho_j) \neq 0$ (for $1 \leq j \leq t$)
by the definition of $\rho_j$.
(Here, $v_{\mathfrak{P}_j}$ is the 
normalized valuation of $K$ with respect to ${\mathfrak{P}_j}$.)
The above equality implies that $a_{r+j}$ must be $0$ for $1 \leq j \leq t$.
Hence we obtain the equalities
\[ \begin{array}{lll}
\beta_1^{a_1} \beta_2^{a_2} \cdots \beta_r^{a_r} = 1
& \mbox{in } \mathcal{U}^1_{\mathfrak{P}_j}, & \mbox{and} \\
\beta_1^{a_1} \beta_2^{a_2} \cdots \beta_r^{a_r} = 1
& \mbox{in } \mathcal{U}^1_{\mathfrak{P}^J_j} & 
\end{array} \]
for all $1 \leq j \leq t$.
Recall that $k^D$ is the decomposition field of $K/\mathbb{Q}$ for $q$.
We denote by $P_1, \ldots, P_u$ the primes of $k^D$ lying above $p$.
As noted before (in the proof for the case NS), we can show
$\mathcal{U}^1_{P_h} \rightarrow 
\prod_{\mathfrak{P} | P_h} \mathcal{U}^1_{\mathfrak{P}}$ is injective.
Hence we see
\[ \beta_1^{a_1} \beta_2^{a_2} \cdots \beta_r^{a_r} = 1
\quad \mbox{in } \mathcal{U}^1_{P_h} \]
for all $1 \leq h \leq u$, and $a_i \neq 0$ with some $i$.
The rest of the proof is quite same as that of for the case NS.
\hfill $\Box$

\par \bigskip

Since $N$ is an abelian extension of $k_\infty$, 
we can take a unique intermediate field $N^+$ of $N/k_\infty$
which satisfies $\Gal (N^+/k_\infty) \cong \mathbb{Z}_p^{\oplus r+s}$.
Similarly, 
we are also able to take a unique intermediate field $N_q^+$ of $N_q/k_\infty$
satisfying $\Gal (N_q^+/k_\infty) \cong \mathbb{Z}_p^{\oplus r}$.
(Note that $N_q^+ \subseteq N^+$, and $N_q^+ = N^+$ for the case NS.)
Then we obtain the following:

\par \bigskip

\noindent
\textbf{Proposition 3.4. }
{\itshape
$N^+/k_\infty$ is unramified outside $\{p,q\}$, 
and a subgroup of $\Gal (N^+/k_\infty)$ generated by the inertia groups 
for the primes lying above $p$ has finite index.
$N_q^+ \cap M_p (k_\infty)/k_\infty$ is a finite extension.
}

\section{Proof of Theorem 1.1}
We will use the same notation and symbols defined in section 3.
Our strategy of the Proof of Theorem 1.1 is similar to that of Theorem A.
However, our situation has a difficulty which comes from the fact that 
$\Gal (M_p (k_\infty)/k_\infty)$ can be non-trivial.
Assume that $K, k$, and $q$ satisfy (A) stated in section 2.

We shall recall and define the following symbols:
\begin{itemize}
\item $M_{p,q} (k_\infty)$ : the maximal pro-$p$ abelian extension of $k_\infty$
unramified outside $\{p,q \}$,
\item $M_p (k_\infty)$ : the maximal pro-$p$ abelian extension of $k_\infty$
unramified outside $p$,
\item $M_q (k_\infty)$ : the maximal pro-$p$ abelian extension of $k_\infty$
unramified outside $q$, 
\item $L (k_\infty)$ : the maximal unramified pro-$p$ abelian extension of $k_\infty$,
\item $\mathfrak{X}_{p,q} (k_\infty) = \Gal (M_{p,q} (k_\infty)/k_\infty)$,
\item $\mathfrak{X}_p (k_\infty) = \Gal (M_p (k_\infty)/k_\infty)$,
\item $X_q (k_\infty) = \Gal (M_q (k_\infty)/k_\infty)$,
\item $X (k_\infty) = \Gal (L (k_\infty)/k_\infty)$.
\end{itemize}

We also define the following notation:
\begin{itemize}
\item $\Gamma = \Gal(K_\infty/K)$ (we often identify $\Gamma$ with $\Gal(k_\infty/k)$.),
\item $\gamma$ : fixed topological generator of $\Gamma$, 
\item $\kappa$ : ($p$-adic) cyclotomic character of $\Gamma$,
\item $\Lambda = \mathbb{Z}_p [[T]] \cong \mathbb{Z}_p [[\Gamma]] : 
1+T \leftrightarrow \gamma$,
\item $\dot{T} = \kappa(\gamma) (1+T)^{-1}-1 \in \Lambda$.
\end{itemize}

We note that $\mathfrak{X}_{p,q} (k_\infty)$, 
$\mathfrak{X}_{p} (k_\infty)$, $X_{q} (k_\infty)$, $X (k_\infty)$, 
and $\Gal (M_{p} (k_\infty)/L(k_\infty))$ are finitely generated 
torsion $\Lambda$-modules.

For a finitely generated torsion $\Lambda$-module $A$, 
we denote by $\charLambda A$ the characteristic ideal of $A$.
For finitely generated torsion $\Lambda$-modules $A$ and $B$,
we write $A \sim B$ when they are pseudo-isomorphic.
We denote by $X (K_\infty)^- := X (K_\infty)^{1-J}$ 
the minus part of $X (K_\infty)$.

We recall the fact that $\mathfrak{X}_p (k_\infty)$ relates to $X (K_\infty)^-$ by 
Kummer duality.
Let $f(T) \in \Lambda$ be a generator of $\charLambda X (K_\infty)^-$.
We note that $f(T)$ is not divisible by $p$ because $K$ is an abelian field
(Ferrero-Washington's theorem \cite{FW}).
It is known that $f(\dot{T}) \in \Lambda$ generates 
$\charLambda \mathfrak{X}_p (k_\infty)$.
By a result of Greenberg \cite{Gre73}, we know that the power of 
$T$ dividing $f(T)$ is $T^s$, 
where $s=0$ for the case NS, and $s=t$ for the case S.
Hence the power of $\dot{T}$ dividing $f(\dot{T})$ is just $\dot{T}^s$.
For the case NS, we see that $f(\dot{T})$ is prime to $\dot{T}$.

For $i=1, \ldots, t$,
let $k_{n,i}$ be the completion of $k_n$ at the unique prime lying above $\mathfrak{p}_i$,
and $\mathcal{U}^1 (k_{n,i})$ the group of principal units in $k_{n,i}$.
Let $\phi(E_{k_n})$ be the diagonal image of $E_{k_n}$ in 
$\prod_i k_{n,i}^{\times}$, and $\mathcal{E}_n$ the 
closure of $\phi(E_{k_n}) \cap \prod_i \mathcal{U}^1 (k_{n,i})$.
We put $\mathcal{U} = \varprojlim \prod_i \mathcal{U}^1 (k_{n,i})$, and 
$\mathcal{E} = \varprojlim \mathcal{E}_n$, where the projective limits are
taken with respect to the norm mappings.
Recall the exact sequence: 
\[ 0 \rightarrow \Gal(M_p (k_\infty)/L(k_\infty)) \rightarrow \mathfrak{X}_p (k_\infty)
\rightarrow X (k_\infty) \rightarrow 0, \]
and the fact that $\Gal(M_p (k_\infty)/L(k_\infty)) \cong \mathcal{U}/\mathcal{E}$
(see \cite[Corollary 13.6]{Was}).
We note that $\varprojlim \mathcal{U}^1 (k_{n,i})$ contains 
$\varprojlim \mu_{p^n} \cong \Lambda/\dot{T}$ for the case S
 (see \cite{Tsuji}, etc.).
Hence $\mathcal{U}$ contains a submodule which is isomorphic to
$(\Lambda/\dot{T})^{\oplus s}$. 
We also note that 
$\mathcal{E}_n$ has no non-trivial $\mathbb{Z}_p$-torsion element 
(see, e.g., \cite[Lemma 3.3]{Tsuji03}).
From the above facts, we obtain the following:

\par \bigskip

\noindent
\textbf{Lemma 4.1. } 
{\itshape
There is a pseudo-isomorphism of finitely generated torsion $\Lambda$-modules: 
\[ \Gal (M_p (k_\infty)/L(k_\infty)) \sim 
(\Lambda/\dot{T})^{\oplus s} \oplus E, \]
where $E$ is an elementary torsion $\Lambda$-module 
(see \cite{Iwa73}, \cite[(5.3.9) Definition]{NSW}, \cite[Chapter 15]{Was})
whose characteristic ideal is prime to $(\dot{T})$.
Moreover, the characteristic ideal of $X(k_\infty)$ is 
prime to $(\dot{T})$.
(See also \cite{Tsuji03}.)
}

\par \bigskip

Let $N^+$ and $N_q^+$ be extensions over $k_\infty$ defined in section 3
(see the paragraph before Lemma 3.4).

\par \bigskip

\noindent
\textbf{Lemma 4.2} (see also Greenberg \cite{Gre77}).
{\itshape
$M_{p,q} (k_\infty) = M_p (k_\infty) N_q^+$.
}

\par \bigskip

\noindent
\textbf{Proof. }
Although this fact is already shown in \cite[pp.148--149]{Gre77}, we will give a detailed proof
for a convenient to the reader
(and our proof is slightly different).
Let $\widetilde{M}_{p,q} (k_\infty)$ be the maximal pro-$p$ extension of $k_\infty$
unramified outside $\{p,q \}$.
By Theorem 3 of Iwasawa \cite{Iwa81} and Ferrero-Washington's theorem \cite{FW}, 
we see that $\Gal (\widetilde{M}_{p,q} (k_\infty)/k_\infty)$ is a free pro-$p$ group
whose minimal number of generators is $\lambda^- + r$,
where $\lambda^- = \zprank X(K_\infty)^-$.
(Note that every prime lying above $q$ actually ramifies in 
$\widetilde{M}_{p,q} (k_\infty)/k_\infty$ by Lemma 3.2.)
By taking the abelian quotient of $\Gal (\widetilde{M}_{p,q} (k_\infty)/k_\infty)$, 
we see that 
$\mathfrak{X}_{p,q} (k_\infty) \cong \mathbb{Z}_p^{\oplus \lambda^-+r}$ 
as a $\mathbb{Z}_p$-module.

On the other hand, we see that $M_p (k_\infty) \cap N_q^+ / k_\infty$ is a finite 
extension by Proposition 3.4.
Hence 
\[ \zprank \Gal(M_p (k_\infty) N_q^+/k_\infty)= \zprank \mathfrak{X}_p (k_\infty)+r
= \lambda^- +r. \]

Since $N_q^+/k_\infty$ is unramified outside $\{p,q\}$, 
we see $M_{p,q} (k_\infty) \supseteq M_p (k_\infty) N_q^+$.
Then we have a surjection of finitely generated $\mathbb{Z}_p$-modules 
$\mathfrak{X}_{p,q} (k_\infty) \rightarrow \Gal(M_p (k_\infty) N_q^+/k_\infty)$
whose kernel is finite.
However $\mathfrak{X}_{p,q} (k_\infty)$ has no non-trivial $\mathbb{Z}_p$-torsion element,
and hence we conclude that $M_{p,q} (k_\infty) = M_p (k_\infty) N_q^+$.
\hfill $\Box$

\par \bigskip

For a finitely generated torsion $\Lambda$-module $A$, we can define the 
``multiplication by $\dot{T}$ endomorphism'' of $A$, 
and we denote by $A[\dot{T}]$ (resp. $A/\dot{T}$) its kernel (resp. cokernel):
\[ 0 \rightarrow A[\dot{T}] \rightarrow A \overset{\dot{T}}{\rightarrow} 
A \rightarrow A/\dot{T} \rightarrow 0. \]

Note that $\Gamma$ acts on $\Gal(M_{p,q} (k_\infty)/L(k_\infty))$ and then 
it is also a finitely generated torsion $\Lambda$-module.
Let $M'$ be the intermediate field of $M_{p,q} (k_\infty)/L(k_\infty)$ corresponding to 
$\dot{T} \Gal(M_{p,q} (k_\infty)/L(k_\infty))$.
Hence $\Gal (M'/L(k_\infty))$ is isomorphic to 
$\Gal(M_{p,q} (k_\infty)/L(k_\infty))/\dot{T}$.

\par \bigskip

\noindent
\textbf{Lemma 4.3. }
{\itshape
$M_q (k_\infty)$ is contained in $M'$.
}

\par \bigskip

\noindent
\textbf{Proof. }
By class field theory, $\Gal(M_q (k_\infty)/L(k_\infty))$ is isomorphic to a 
quotient of $\varprojlim \widehat{(O_{k_n}/q)^{\times}}$.
As a $\Lambda$-module, $\varprojlim \widehat{(O_{k_n}/q)^{\times}}$
is isomorphic to $(\Lambda/ \dot{T})^{\oplus r}$.
Hence $\dot{T}$ annihilates $\Gal(M_q (k_\infty)/L(k_\infty))$, 
and then $\Gal(M_q (k_\infty)/L(k_\infty))/\dot{T} =
\Gal(M_q (k_\infty)/L(k_\infty))$.
From the restriction map 
\[ \Gal(M_{p,q} (k_\infty)/L(k_\infty)) \rightarrow 
\Gal(M_q (k_\infty)/L(k_\infty)) \rightarrow 0, \]
we obtain a surjection
\[ \Gal(M_{p,q} (k_\infty)/L(k_\infty))/\dot{T} \rightarrow 
\Gal(M_q (k_\infty)/L(k_\infty)) \rightarrow 0. \]
By the definition of $M'$, we see $M_q (k_\infty) \subseteq M'$.
\hfill $\Box$

\par \bigskip

\noindent
\textbf{Lemma 4.4. }
{\itshape
$L(k_\infty) N^+$ is contained in $M'$.
}

\par \bigskip

\noindent
\textbf{Proof. }
Note that $\Gal(L(k_\infty) N^+/L(k_\infty)) \cong 
\Gal(N^+/N^+ \cap L(k_\infty))$, and $\Gal(N^+/N^+ \cap L(k_\infty))$ is 
a subgroup of $\Gal(N^+/k_\infty)$.
By the construction of $N^+$, we see that $\dot{T}$ annihilates 
$\Gal(N^+/k_\infty)$, and hence it also annihilates 
$\Gal(L(k_\infty)N^+/L(k_\infty))$.
The rest of the proof is similar to that of Lemma 4.3.
\hfill $\Box$

\par \bigskip

\noindent
\textbf{Lemma 4.5. }
{\itshape
$M'/L(k_\infty)N^+$ is a finite extension.
}

\par \bigskip

\noindent
\textbf{Proof. }
We shall show that 
$\zprank \Gal(M'/L(k_\infty)) = \zprank \Gal(L(k_\infty)N^+/L(k_\infty))$.
By Proposition 3.3, we see that $N^+ \cap L(k_\infty)/k_\infty$ is a finite extension.
Hence $\zprank \Gal(L(k_\infty)N^+/L(k_\infty))$ is equal to
$\zprank \Gal(N^+/k_\infty)=r+s$.

On the other hand, $M_{p,q} (k_\infty) = M_p (k_\infty) N_q^+$ by Lemma 4.2, 
and $M_p (k_\infty) \cap N_q^+/k_\infty$ is a finite extension by Proposition 3.4.
By using Lemma 4.1, 
we can obtain the following pseudo-isomorphisms:
\[ \mathfrak{X}_{p,q} (k_\infty) \sim 
\mathfrak{X}_p (k_\infty) \oplus \Gal(N_q^+/k_\infty) 
\sim (\Lambda/\dot{T})^{\oplus r+s} \oplus E', \]
where $E'$ is an elementary torsion $\Lambda$-module 
whose characteristic ideal is prime to $(\dot{T})$.
Hence, $\zprank \mathfrak{X}_{p,q} (k_\infty)/\dot{T}=r+s$.

The following exact sequence: 
\[ 0 \rightarrow \Gal(M_{p,q} (k_\infty) /L(k_\infty)) \rightarrow 
\mathfrak{X}_{p,q} (k_\infty) \rightarrow X (k_\infty) \rightarrow 0 \]
induces the exact sequence: 
\[ X(k_\infty) [\dot{T}] \rightarrow 
\Gal(M_{p,q} (k_\infty) /L(k_\infty))/\dot{T} \rightarrow 
\mathfrak{X}_{p,q} (k_\infty)/\dot{T} \rightarrow X (k_\infty)/\dot{T} 
\rightarrow 0. \]
Since $\charLambda X(k_\infty)$ is prime to $(\dot{T})$ by Lemma 4.1,
both of $X(k_\infty) [\dot{T}]$ and $X (k_\infty)/\dot{T}$ are finite.
Hence $\zprank \Gal(L(k_\infty)N^+/L(k_\infty))$ is equal to 
$\zprank \mathfrak{X}_{p,q} (k_\infty)/\dot{T} = \zprank \Gal(M'/L(k_\infty))$.
\hfill $\Box$

\par \bigskip

For a Galois group $G$ appeared below, 
we denote $\mathcal{I} (G)$ by the subgroup of $G$ generated by the 
inertia groups for all primes lying above $p$.

\par \bigskip

\noindent
\textbf{Lemma 4.6. }
{\itshape
$\zprank \mathcal{I} (\Gal(N^+L(k_\infty)/L(k_\infty))) = r+s$.
}

\par \bigskip

\noindent
\textbf{Proof. }
We shall take 
a prime $\mathcal{P}$ of $k_\infty$ lying above $p$.
Let $I_{\mathcal{P}}$ be the inertia subgroup of $\Gal(N^+/k_\infty)$ for $\mathcal{P}$.
Similarly, let $I'_{\mathcal{P}}$ be the inertia subgroup of 
$\Gal(N^+L(k_\infty)/k_\infty)$ for $\mathcal{P}$.
Then the restriction map induces a surjection
$I'_{\mathcal{P}} \to I_{\mathcal{P}}$.
Hence there is a surjection
$\mathcal{I} (\Gal(N^+L(k_\infty)/k_\infty)) \to \mathcal{I} (\Gal(N^+/k_\infty))$.
By Proposition 3.4, $\zprank \mathcal{I} (\Gal(N^+/k_\infty))=r+s$.
We see that $\zprank \mathcal{I} (\Gal(N^+L(k_\infty)/k_\infty)) \geq r+s$.
Since $L(k_\infty)/k_\infty$ is an unramified extension, 
$\mathcal{I} (\Gal(N^+L(k_\infty)/k_\infty))$ is contained in 
$\Gal(N^+L(k_\infty)/L(k_\infty))$.
By these results, 
we see that $\zprank \mathcal{I} (\Gal(N^+L(k_\infty)/L(k_\infty))) = r+s$.
\hfill $\Box$

\par \bigskip

We shall finish to prove Theorem 1.1.
By Lemma 4.6, $\zprank \mathcal{I} (\Gal(N^+L(k_\infty)/L(k_\infty))) = r+s$.
Moreover, $\zprank \mathcal{I} (\Gal(M'/L(k_\infty)))$ 
is also $r+s$ because $M' / N^+ L(k_\infty)$ is a finite extension
(Lemma 4.5).
From the proof of Lemma 4.5, we see that $\zprank \Gal (M'/L(k_\infty))$
is $r+s$.
Then $\mathcal{I} (\Gal(M'/L(k_\infty)))$ is a finite index subgroup 
of $\Gal (M'/L(k_\infty))$.
By Lemma 4.3, $M_q (k_\infty)$ is an intermediate field of $M'/L(k_\infty)$.
Since $M_q (k_\infty)/L(k_\infty)$ is unramified at all primes lying above $p$,
we can see that $M_q (k_\infty)$ is contained in the fixed field 
of $\mathcal{I} (\Gal(M'/L(k_\infty)))$.
This implies that $M_q (k_\infty)/L(k_\infty)$ is a finite extension.

We have shown Theorem 1.1 for $k$ and $q$ satisfying (A).
Then, as noted in section 2, we obtain Theorem 1.1 for general $k$ and $q$.
\hfill $\Box$

\section{Slight generalization of Theorem 1.1}
In this section, we shall give a simple remark that the converse of 
Lemma 2.1 holds under a (strict) condition.
(In general, the converse of Lemma 2.1 does not hold.)

\par \bigskip

\noindent
\textbf{Lemma 5.1. }
{\itshape
Let $k'/k$ be a finite extension of algebraic number fields 
satisfying $k' \cap k_\infty = k$.
Let 
\[ I_n : \widehat{(O_{k_n}/q)^{\times}} \rightarrow \widehat{(O_{k'_n}/q)^{\times}} \]
be the homomorphism induced from the natural embedding.
If $\Gal(M_q (k_\infty) / L(k_\infty))$ is finite 
and $I_n$ is an isomorphism for all $n$, 
then $\Gal(M_q (k'_\infty) / L(k'_\infty))$ is also finite.
}

\par \bigskip

\noindent
\textbf{Proof. }
By the assumption, we obtain the following isomorphism
\[ I : \varprojlim \widehat{(O_{k_n}/q)^{\times}} \rightarrow 
\varprojlim \widehat{(O_{k'_n}/q)^{\times}}. \]
Hence the homomorphism 
\[ \Gal(M_q (k_{\infty}) / L(k_{\infty})) \rightarrow 
\Gal(M_q (k'_{\infty}) / L(k'_{\infty})) \]
induced from $I$ is surjective.
The assertion follows.
\hfill $\Box$

\par \bigskip

We shall give an example satisfying the assumption that $I_n$ is an 
isomorphism for all $n$.
Let $k$ be an algebraic number field, and $k'$ a quadratic extension of $k$.
Assume that every prime of $k$ lying above $q$ is inert in $k'_{\infty}$
(that is, every prime lying above $q$ is inert in $k'$ and $k_\infty$).
Let $\mathfrak{q}_1, \ldots, \mathfrak{q}_r$ be the primes of $k$ lying 
above $q$.
We also assume that $p$ divides $N(\mathfrak{q}_i)-1$ for all $i$, 
where $N(\mathfrak{q}_i)$ is the absolute norm of $\mathfrak{q}_i$.
Then $\widehat{(O_{k_n}/\mathfrak{q}_i)^{\times}}$ is not trivial 
for all $i, n$.
Under these assumptions, we obtain that 
\[ N(\mathfrak{q}_i O_{k'_n})-1 = N(\mathfrak{q}_i O_{k_n})^2-1 =
(N(\mathfrak{q}_i O_{k_n})-1)(N(\mathfrak{q}_i O_{k_n})+1). \]
We assumed that $p$ is odd, and hence $N(\mathfrak{q}_i O_{k_n})+1$ 
is prime to $p$.
This implies that $|\widehat{(O_{k'_n}/\mathfrak{q}_i)^{\times}}|
= |\widehat{(O_{k_n}/\mathfrak{q}_i)^{\times}}|$ for all $i, n$.
From this, we obtain that 
\[ |\widehat{(O_{k'_n}/q)^{\times}}| = \prod_{i=1}^r 
|\widehat{(O_{k'_n}/\mathfrak{q}_i)^{\times}}| = 
\prod_{i=1}^r 
|\widehat{(O_{k_n}/\mathfrak{q}_i)^{\times}}| =
|\widehat{(O_{k_n}/q)^{\times}}|. \]
Since $I_n$ is injective, we see that $I_n$ is an isomorphism for all $n$.
Under the above assumptions, if $\Gal(M_q (k_\infty) / L(k_\infty))$ 
is finite, then $\Gal(M_q (k'_\infty) / L(k'_\infty))$ is also finite
by Lemma 5.1.
(For the case that $k'/k$ is an 
imaginary quadratic extension of $\mathbb{Q}$, 
see also \cite{Salle}, \cite{IMO}.)

The above lemma implies that Theorem 1.1 can be generalized for some 
non-abelian fields.

\section{$\mathbb{Z}_p$-rank of $S$-ramified Iwasawa modules}
We shall lead a formula of the $\mathbb{Z}_p$-rank of $S$-ramified 
Iwasawa modules (for general $S$) from Theorem 1.1.
As same as Theorem 1.1,
the strategy of our proof is quite similar to that of given in \cite{IMO}.

In this section, we will use the following notation 
(similar to Greither's \cite{Greither} or Tsuji's \cite{Tsuji} 
but slightly different):
\begin{itemize}
\item $p$ : fixed odd prime,
\item $S$ : finite set of \textit{rational} primes which does not include $p$,
\item $F$ : finite abelian extension of $\mathbb{Q}$ unramified at $p$,
\item $K = F (\mu_p)$,
\item $K_n = K (\mu_{p^{n+1}})$,
\item $K_\infty= \cup_{n \geq 1} K_n$ : the cyclotomic $\mathbb{Z}_p$-extension of $K$,
\item $G = \Gal (K_\infty/\mathbb{Q}_\infty) \cong \Gal (K/\mathbb{Q})$,
\item $\Gamma = \Gal(K_\infty /K)$,
\item $G_p$ : Sylow $p$-subgroup of $G$,
\item $G_0$ : non-$p$-part of $G$
(the maximal subgroup of $G$ consists of the elements having prime to $p$ order),
\item $\gamma$ : fixed topological generator of $\Gamma$,
\item $\kappa$ : ($p$-adic) cyclotomic character,
\item $\omega$ : ($p$-adic) Teichm{\"u}ller character,
\item $J \in G$ : complex conjugation.
\end{itemize}

Let $\chi$ be a $p$-adic character of $G$.
We denote by 
$\mathbb{Q}_p (\chi)$ the extension of $\mathbb{Q}_p$ 
by adjoining the values of $\chi$,
and $O_\chi$ the valuation ring of $\mathbb{Q}_p (\chi)$.
We put $d_\chi = [\mathbb{Q}_p (\chi) : \mathbb{Q}_p ]$.
Let $\underline{O_\chi}$ be a free rank one $O_\chi$-module
such that $\sigma \in G$ acts as $\chi (\sigma)$.
For a $\mathbb{Z}_p [G]$-module $M$, we put
$M_{\chi} = M \otimes_{\mathbb{Z}_p [G]} \underline{O_\chi}$,
which is called the ``$\chi$-quotient'' in \cite{Tsuji}
(or the ``$\chi$-part'' in \cite{Greither}). 
The functor taking the $\chi$-quotient is right exact.
We also put 
\[ e_{\chi} = \frac{1}{|G|} \sum_{\sigma \in G} 
\mathrm{tr}_{\mathbb{Q}_p (\chi)/\mathbb{Q}_p} (\chi(\sigma)) \sigma^{-1} 
\in \mathbb{Q}_p [G]. \]
If $p$ does not divide $|G|$, then $M_\chi \cong e_\chi M$.
In general, we see 
\[ M_\chi \otimes_{\mathbb{Z}_p} \mathbb{Q}_p
\cong (M \otimes_{\mathbb{Z}_p} \mathbb{Q}_p)_\chi 
\cong e_\chi (M \otimes_{\mathbb{Z}_p} \mathbb{Q}_p). \]
For more informations about the $\chi$-quotient, see \cite{Greither}, \cite{Tsuji}
for example.

We also give some simple remarks.
For a $\mathbb{Z}_p [G]$-module $M$,
we put $M^{\pm}=M^{1 \pm J}$.
Since $p$ is odd, we have a decomposition $M \cong M^+ \oplus M^-$.
For a character $\chi$ of $G$, we see that 
\[ \begin{array}{rcl}
M_\chi & \cong & (M^+ \oplus M^-)_\chi \\
& = & (M^+ \oplus M^-) \otimes_{\mathbb{Z}_p [G]} \underline{O_\chi} \\
& \cong & (M^+ \otimes_{\mathbb{Z}_p [G]} \underline{O_\chi})
\oplus (M^- \otimes_{\mathbb{Z}_p [G]} \underline{O_\chi}) \\
& = & M^+_\chi \oplus M^-_\chi. 
\end{array} \]
We claim that if $\chi$ is odd, then $M^+_\chi$ is trivial.
Let 
\[ (a \otimes b) \in M^+ \otimes_{\mathbb{Z}_p [G]}
\underline{O_\chi} = M^+_\chi. \]
Note that $J$ acts trivially on $M^+$ and acts as $-1$ on $\underline{O_\chi}$.
Hence the equality
\[ (a \otimes b) = (Ja \otimes b) = (a \otimes Jb) = (a \otimes -b) \]
implies that $(a \otimes 2b) = 2(a \otimes b)=0$.
Since $p$ is odd, we obtain the claim.
Similarly, we can see that if $\chi$ is even, then 
$M^-_\chi$ is trivial.

For a rational prime $q$ distinct from $p$, we put
\[ R_q = \varprojlim \widehat{(O_{K_n} /q)^{\times}}. \]
Let $r$ be the number of primes of $K_\infty$ lying above $q$.
Then $\zprank R_q = r$.
Since $q$ is a rational prime, $G$ acts on $R_q$.
We shall determine the $\mathbb{Z}_p$-rank of $(R_q)_\chi$.

First, we assume that $q$ is unramified in $K$
(i.e., the conductor of $F$ is prime to $q$).
Let $D$ be the decomposition subgroup of $\Gal(K_\infty/\mathbb{Q})$ for $q$.
Then we can write $D \cong D_p \times D_0$,
where $D_p \cong \mathbb{Z}_p$ and
$D_0$ is a finite cyclic group whose order is prime to $p$. 
We may regard $D_0$ as a subgroup of $G_0$.
Note that $\Gal (K_\infty/\mathbb{Q})$ is isomorphic to $\Gamma \times G_p \times G_0$.
Then we can take a generator of $D_p$ of the from $\gamma^{p^m} \sigma_p$ with 
some $m \geq 0$ and $\sigma_p \in G_p$.
We also take a generator $\sigma_0 \in G_0$ of $D_0$.
Hence $D$ is a procyclic group generated by $\gamma^{p^m} \sigma_p \sigma_0$.

In the above choice of the generator of $D_p$, we can see that 
$m$ and $\sigma_p$ is uniquely determined.
(Since $D_p \cong \mathbb{Z}_p$, every generator of $D_p$ is written by the
from $(\gamma^{p^m} \sigma_p)^\alpha$ with $\alpha \in \mathbb{Z}_p^{\times}$.)

\par \bigskip

\noindent
\textbf{Lemma 6.1. }
{\itshape
$R_q$ is a cyclic $\mathbb{Z}_p [G] [[\Gamma ]]$-module.
}

\par \bigskip

\noindent
\textbf{Proof. }
Fix a sufficiently large integer $n_0$ such that every prime in $K_{n_0}$ lying above $q$ 
remains prime in $K_m$ for all $m \geq n_0$.
Let $n$ be an integer which satisfies $n \geq n_0$.
We put $G^{(n)} = \Gal (K_n/\mathbb{Q})$.
Let $\mathfrak{q}$ be a prime in $K_n$ lying above $q$.
We remark that $\mu_{p^{n+1}} \subset K_n$ and $\mu_{p^{n+2}} \not\subset K_n$.
Under the assumption on $n$, we can see that the Sylow $p$-subgroup of 
$(O_{K_n}/\mathfrak{q})^{\times}$ is
generated by $\zeta_{p^{n+1}} \pmod{\mathfrak{q}}$ with a generator 
$\zeta_{p^{n+1}}$ of $\mu_{p^{n+1}}$.
Let $\{ \mathfrak{q}_1, \mathfrak{q}_2, \ldots, \mathfrak{q}_r \}$ be the set of 
primes of $K_n$ lying above $q$.
We assumed that $q$ is unramified in $K_\infty/\mathbb{Q}$, then 
$q O_{K_n} = \mathfrak{q}_1 \mathfrak{q}_2 \cdots \mathfrak{q}_r$.
We note that the action of $G^{(n)}$ on 
$\{ \mathfrak{q}_1, \mathfrak{q}_2, \ldots, \mathfrak{q}_r \}$ 
is transitive.
Take an element $\alpha_n$ of $O_{K_n}$ which satisfies
\[ \alpha_n \equiv \zeta_{p^{n+1}} \pmod{\mathfrak{q}_1}, \quad
\alpha_n \equiv 1 \pmod{\mathfrak{q}_2}, \quad
\ldots, \quad \alpha_n \equiv 1 \pmod{\mathfrak{q}_r}. \]
Then $\alpha_n \pmod{q}$ is a generator of the Sylow $p$-subgroup of $(O_{K_n}/q)^{\times}$ 
as a $\mathbb{Z}_p [G^{(n)}]$-module.
Hence the Sylow $p$-subgroup of $(O_{K_n}/q)^{\times}$ (which is isomorphic to 
$\widehat{(O_{K_n}/q)^{\times}}$) is a cyclic $\mathbb{Z}_p [G^{(n)}]$-module.

We can choose a suitable set of generators $\{ \alpha_n \}$ 
such that $N_{K_m/K_n} (\alpha_m) \equiv \alpha_n \pmod{q}$ 
for all $m>n \geq n_0$.
Hence we obtain the following commutative diagram with exact raws:
\[ \begin{array}{ccccc}
\mathbb{Z}_p [G^{(m)}] & \rightarrow & \widehat{(O_{K_m}/q)^{\times}} & \rightarrow & 0 \\
\downarrow & & \downarrow & &  \\
\mathbb{Z}_p [G^{(n)}] & \rightarrow & \widehat{(O_{K_n}/q)^{\times}} & \rightarrow & 0,
\end{array} \]
where the left vertical mapping is induced from the restriction mapping,
and the right vertical mapping is induced from the norm mapping.
Since $\varprojlim \mathbb{Z}_p [G^{(n)}] \cong 
\mathbb{Z}_p [G] [[\Gamma]]$, 
we obtain the assertion.
\hfill $\Box$

\par \bigskip

Hence there is a surjection $\varphi : \mathbb{Z}_p [G] [[\Gamma ]] \to R_q$.
We note that $\gamma^{p^m} \sigma_p \sigma_0$ acts on $R_q$ as 
$\kappa (\gamma^{p^m} \sigma_p \sigma_0)$.
(Recall that $\kappa$ is the cyclotomic character.)
Then the kernel of $\varphi$ contains an ideal generated by 
$\gamma^{p^m} \sigma_p \sigma_0 - \kappa (\gamma^{p^m} \sigma_p \sigma_0)$.

By taking the $\chi$-quotient, we obtain a surjection 
$\varphi_\chi : O_\chi [[\Gamma ]] \to (R_q)_\chi$, and 
the kernel of $\varphi_\chi$ contains 
$\chi(\sigma_p \sigma_0) \gamma^{p^m} - \kappa (\gamma^{p^m} \sigma_p \sigma_0)$.
We may regard $(R_q)_\chi$ as a $O_\chi [[T]]$-module
via the isomorphism $O_\chi [[\Gamma ]] \cong O_\chi [[T]]$
with $\gamma \mapsto 1+T$.
We put $\Lambda_\chi =O_\chi [[T]]$, 
and $\kappa_0 = \kappa (\gamma) \in 1+p \mathbb{Z}_p$.
Then we see that $(R_q)_\chi$ is annihilated by 
\[ f_{q,\chi} (T) = (1+T)^{p^m} - \chi^{-1} (\sigma_p \sigma_0) 
\kappa (\sigma_p \sigma_0) \kappa_0^{p^m}
\in \Lambda_\chi. \]

Let $\mathfrak{P}$ be the maximal ideal of $O_\chi$.
Since $\kappa_0 \in 1+p \mathbb{Z}_p$, if 
\[ \chi^{-1} (\sigma_p \sigma_0) \kappa (\sigma_p \sigma_0) \not\equiv 
1 \pmod{\mathfrak{P}}, \]
then $f_{q,\chi} (T)$ is a unit polynomial,
and hence $(R_q)_\chi$ is trivial.
We see
\[ \chi^{-1} (\sigma_p \sigma_0) \kappa (\sigma_p \sigma_0)
= \chi^{-1} \kappa (\sigma_0) \chi^{-1} \kappa (\sigma_p) 
= \chi^{-1} \omega (\sigma_0) \chi^{-1} \kappa (\sigma_p). \]
Note that $\chi^{-1} \kappa (\sigma_p)$ is a $p$-power root of unity, and then 
it is congruent to $1$ modulo $\mathfrak{P}$.
Moreover, $\chi^{-1} \omega (\sigma_0)$ is a root of unity whose order is prime to $p$.
Then $\chi^{-1} \omega (\sigma_0) \equiv 1 \pmod{\mathfrak{P}}$ 
if and only if $\chi^{-1} \omega (\sigma_0)=1$.
Consequently, we showed that if $\chi^{-1} \omega (\sigma_0) \neq 1$, 
then $(R_q)_\chi$ is trivial.

We can see that the number of characters $\chi$ satisfying 
$\chi^{-1} \omega (\sigma_0) = 1$ is just $|G/ D_0|$.
(It is equal to the number of characters $\chi'$ of $G$ 
satisfying $\chi' (D_0)=1$.)
Though if $(R_q)_\chi$ is non-trivial, it is annihilated by $f_{q,\chi} (T)$,
and then $\zprank (R_q)_\chi \leq d_\chi p^m$.
By considering these facts, we obtain the inequality:
\[ r = \zprank R_q = \sum_\chi \zprank (R_q)_\chi
\leq \sum_\chi d_\chi p^m= |G/D_0| \times p^m, \]
where $\chi$ runs all representatives of the conjugacy classes
satisfying $\chi^{-1} \omega (\sigma_0) = 1$
in the above sums.
(We give some remarks. The second equation follows from the 
fact that $R_q \otimes_{\mathbb{Z}_p} \mathbb{Q}_p
\cong \bigoplus_\chi (R_q)_\chi \otimes_{\mathbb{Z}_p} \mathbb{Q}_p$.
Moreover, $\sum_\chi d_\chi = |G/D_0|$, and
$m$ is independent of $\chi$).

We claim that $r = |G/D_0| \times p^m$.
Let $K^D$ be the decomposition field of $K_\infty/\mathbb{Q}$ for $q$.
Then the $p$-part of $[K^D : \mathbb{Q}]$ is equal to the $p$-part of
$[K_m : \mathbb{Q} ]$.
Hence this is equal to $|G_p| \times p^m$.
On the other hand, the non-$p$-part of $[K^D : \mathbb{Q}]$ is equal to 
$|G_0/D_0|$.
We see 
\[ [K^D : \mathbb{Q}] = |G_p| \times p^m \times |G_0/D_0|
= |G/D_0| \times p^m. \]
Since $r = [K^D : \mathbb{Q}]$, the claim follows.

From this claim, we see that the above inequality is just an equality.
Hence for all character $\chi$ satisfying $\chi^{-1} \omega (\sigma_0) = 1$, 
the $\mathbb{Z}_p$-rank of $(R_q)_\chi$ is $d_\chi p^m$.
We also note that $\kappa (\sigma_p)=1$
because $\sigma_p$ fixes all elements of $\mu_{p^n}$ for all $n$. 
Hence, when $\chi^{-1} \omega (\sigma_0) = 1$, 
we can write 
\[ f_{q,\chi} (T) = (1+T)^{p^m} - \chi^{-1} (\sigma_p) \kappa_0^{p^m}. \]

Next, we consider the case that $q$ is ramified in $K$.
Let $I$ be the inertia subgroup of $\Gal (K/\mathbb{Q})$ for $q$,
and $K^I$ the inertia field of $K/\mathbb{Q}$ for $q$.
We remark that all primes lying above $q$ are totally 
ramified in $K_n/ K^I_n$.
Hence $\widehat{(O_{K_n} /q)^{\times}} \cong \widehat{(O_{K^I_n} /q)^{\times}}$
for all $n$.
We put
\[ R^I_q = \varprojlim \widehat{(O_{K^I_n} /q)^{\times}}. \]
Since $q$ is unramified in $\mathbb{Q} (\mu_{p^\infty})$ 
(where $\mu_{p^\infty} = \bigcup_{n \geq 1} \mu_{p^n}$),
we see that $K_\infty^I$ contains $\mu_{p^\infty}$.
Then the proof of Lemma 6.1 also works for $K_\infty^I$.

Let $\chi$ be a character of $G$.
If $\chi (I)=1$, then $\chi$ is also a character of 
$\Gal(K^I / \mathbb{Q})$, and hence 
$(R_q)_\chi \cong (R^I_q)_\chi$.
From this, if $\chi (I) \neq 1$, we see that 
$(R_q)_\chi$ is finite because
\[ \zprank R_q = \zprank R^I_q 
= \sum_{\chi(I)=1} \zprank (R^I_q)_\chi 
= \sum_{\chi(I)=1} \zprank (R_q)_\chi.\]

We also determine the structure of $(R_q)_\chi$ in general case.
Assume that $\chi(I)=1$.
Then $\chi (\sigma)$ for $\sigma \in \Gal(K^I/\mathbb{Q}) \cong G/I$ is 
well defined.
Repeating the argument given in the unramified case for $K^I$,
we can take $\sigma_0$ and $\sigma_p$ for $q$.
(They are determined modulo $I$, 
and $\sigma_p \pmod{I}$ is uniquely determined.
Hence $\chi(\sigma_p)$ is dependent only on $q$.)
We also assume that $\chi^{-1} \omega (\sigma_0)=1$.
Since $\mathbb{Z}_p [G/I]_\chi \cong O_\chi$, 
we can take 
\[ f_{q,\chi} (T) = (1+T)^{p^m} - \chi^{-1} (\sigma_p) \kappa_0^{p^m} \]
as an element of $\Lambda_\chi =O_\chi [[T]]$.
We see that $(R_q)_\chi$ is annihilated by $f_{q,\chi} (T)$ 
because $(R_q)_\chi \cong (R^I_q)_\chi$.

As a consequence, we obtained the following:

\par \bigskip

\noindent
\textbf{Proposition 6.2. }
{\itshape
Let $\chi$ be a character of $G$.
Then $(R_q)_\chi  \otimes_{\mathbb{Z}_p} \mathbb{Q}_p$ is non-trivial if and only if
$\chi$ satisfies $\chi(I)=1$ and $\chi^{-1} \omega (\sigma_0)=1$.
Moreover, if $(R_q)_\chi \otimes_{\mathbb{Z}_p} \mathbb{Q}_p$ is non-trivial, then 
\[ (R_q)_\chi \otimes_{\mathbb{Z}_p} \mathbb{Q}_p \cong 
\Lambda_\chi/ f_{q,\chi} (T) \otimes_{\mathbb{Z}_p} \mathbb{Q}_p, \]
and $\zprank (R_q)_\chi = d_\chi p^m$.
}

\par \bigskip

By class field theory, we have the following exact sequence: 
\[ E_\infty \rightarrow R_q \rightarrow 
\Gal (M_q (K_\infty)/L(K_\infty)) \rightarrow 0, \]
where $E_\infty = \varprojlim \widehat{E_{K_n}}$.
Assume that $\chi$ is a non-trivial \textit{even} character of $G$ 
satisfying $\chi(I)=1$ and $\chi^{-1} \omega (\sigma_0)=1$.
By taking the $\chi$-quotient (it is right exact), 
we see
\[ (E_\infty)_\chi \rightarrow (R_q)_\chi \rightarrow 
\Gal (M_q (K_\infty)/L(K_\infty))_\chi \rightarrow 0 \]
is exact.
Since $(R_q)_\chi$ is annihilated by $f_{q,\chi} (T)$, 
we obtain the exact sequence: 
\[ (E_\infty)_\chi/f_{q,\chi} (T) \rightarrow (R_q)_\chi 
\rightarrow \Gal (M_q (K_\infty)/L(K_\infty))_\chi \rightarrow 0. \]
By tensoring with $\mathbb{Q}_p$, we also obtain the exact sequence: 
\[ (E_\infty)_\chi/f_{q,\chi} (T) \otimes_{\mathbb{Z}_p} \mathbb{Q}_p
\rightarrow (R_q)_\chi \otimes_{\mathbb{Z}_p} \mathbb{Q}_p
\rightarrow \Gal (M_q (K_\infty)/L(K_\infty))_\chi \otimes_{\mathbb{Z}_p} \mathbb{Q}_p
\rightarrow 0. \]

\par \bigskip

\noindent
\textbf{Proposition 6.3. }
{\itshape
For every non-trivial even character $\chi$ of $G$
satisfying $\chi (I)=1$ and $\chi^{-1} \omega (\sigma_0)=1$,
the mapping 
\[ (E_\infty)_\chi /f_{q,\chi} (T) \otimes_{\mathbb{Z}_p} \mathbb{Q}_p
\rightarrow (R_q)_\chi \otimes_{\mathbb{Z}_p} \mathbb{Q}_p \]
appeared above is an isomorphism.
}

\par \bigskip

\noindent
\textbf{Proof. }
As we noted before, we have a decomposition
\[ \Gal (M_q (K_\infty)/L(K_\infty))_\chi \cong 
\Gal (M_q (K_\infty)/L(K_\infty))^+_\chi \oplus \Gal (M_q (K_\infty)/L(K_\infty))^-_\chi. \]
Moreover, we already know that 
$\Gal (M_q (K_\infty)/L(K_\infty))^-_\chi$ is trivial for every even character $\chi$.
Let $K^+$ be the maximal real subfield of $K$.
Since $p$ is odd, we see that $\Gal (M_q (K_\infty)/L(K_\infty))^+$ is isomorphic to
$\Gal (M_q (K^+_\infty)/L(K^+_\infty))$.
Hence we obtain 
\[ \Gal (M_q (K_\infty)/L(K_\infty))_\chi \cong 
\Gal (M_q (K^+_\infty)/L(K^+_\infty))_\chi \]
for every $\chi$ satisfying the assumption.
(We note that $\chi$ can be viewed as a character of $\Gal(K^+_\infty/\mathbb{Q}_\infty)$.)

By Theorem 1.1, we see that $\Gal (M_q (K^+_\infty)/L(K^+_\infty))$ is finite,
and hence we see that 
$\Gal (M_q (K_\infty)/L(K_\infty))_\chi \otimes_{\mathbb{Z}_p} \mathbb{Q}_p$ is trivial.
From this, we have a surjection
\[ (E_\infty)_\chi /f_{q,\chi} (T) \otimes_{\mathbb{Z}_p} \mathbb{Q}_p
\rightarrow (R_q)_\chi \otimes_{\mathbb{Z}_p} \mathbb{Q}_p. \]
By Proposition 6.2, $\dim_{\mathbb{Q}_p}  
(R_q)_\chi \otimes_{\mathbb{Z}_p} \mathbb{Q}_p = d_\chi p^m$.
We shall calculate the dimension of 
$(E_\infty)_\chi/f_{q,\chi} (T) \otimes_{\mathbb{Z}_p} \mathbb{Q}_p$.

Let $\mathcal{U}$ be the projective limit of semi local units in $K_\infty/K$ 
for the primes lying above $p$, and
$\mathcal{E}$ the closure of the diagonal image of global units.
(For the precise definition, see section 4.)
Since Leopoldt's conjecture is valid for all $K_n$ (see \cite{Brumer}, \cite{Was}), 
we see that $\mathcal{E}$ is
isomorphic to $E_\infty$.
It is known that $\mathcal{E}_\chi \otimes_{\mathbb{Z}_p} \mathbb{Q}_p$ is a 
free cyclic $\Lambda_\chi \otimes_{\mathbb{Z}_p} \mathbb{Q}_p$-module.
(See \cite[Lemma 3.5]{Tsuji03}, \cite{Greither}.)
Then we see 
\[ \begin{array}{rcl}
(\mathcal{E}_\chi/f_{q,\chi} (T)) \otimes_{\mathbb{Z}_p} \mathbb{Q}_p
 & \cong & (\mathcal{E}_\chi \otimes_{\mathbb{Z}_p} \mathbb{Q}_p)/
(f_{q,\chi} (T) \otimes 1) \\
 & \cong & (\Lambda_\chi \otimes_{\mathbb{Z}_p} \mathbb{Q}_p)/
(f_{q,\chi} (T) \otimes 1) \\
& \cong & 
(\Lambda_\chi/f_{q,\chi} (T)) \otimes_{\mathbb{Z}_p} \mathbb{Q}_p. 
\end{array} \]
Hence we showed that $\dim_{\mathbb{Q}_p} 
(E_\infty)_\chi/f_{q,\chi} (T) \otimes_{\mathbb{Z}_p} \mathbb{Q}_p
= d_\chi p^m$.
This implies the assertion.
\hfill $\Box$

\par \bigskip

Here we shall state our main result.
Let $\chi$ be an arbitrary character of $G$.
Let $S$ be a finite set of rational primes which does not include $p$. 
We put $X_S (K_\infty)=\Gal (M_S (K_\infty)/K_\infty)$, where 
$M_S (K_\infty)/K_\infty$ is the maximal abelian pro-$p$ extension
unramified outside $S$.
In this case, $G$ acts on $X_S (K_\infty)$ and hence its $\chi$-quotient 
can be considered.
We shall give a formula of $\zprank X_S (K_\infty)_\chi$.
For a prime $q \in S$, let $I_q$ be the inertia subgroup of $G$ for $q$.
We also write $\sigma_{p,q}, \sigma_{0,q}, m_q$
as $\sigma_p, \sigma_0, m$ for $q$ (defined before), respectively.
(Recall that $\sigma_{p,q}$ and $\sigma_{0,q}$ are determined modulo $I_q$.)
We also recall that $\kappa_0 = \kappa (\gamma)$.

\par \bigskip

\noindent
\textbf{Theorem 6.4. }
{\itshape
We put
\[ S_\chi = \{ q \in S \; | \; \chi (I_q) =1, \; \chi^{-1} \omega (\sigma_{0,q})=1 \}, \]
\[ f_{q,\chi} (T) = (1+T)^{p^{m_q}} - \chi^{-1} (\sigma_{p,q}) \kappa_0^{p^{m_q}}
\in O_\chi [[T]], \]
and $F(T)=\mathrm{lcm}_{q \in S_{\chi}} f_{q, \chi} (T)$.
If $S_\chi$ is not empty, then 
\[ \zprank X_S (K_\infty)_\chi = \zprank X (K_\infty)_\chi
+ \sum_{q \in S_{\chi}} d_\chi p^{m_q} - P_\chi, \]
where 
\[ P_\chi = \left\{ \begin{array}{ll}
1 & \mbox{($\chi= \omega$)} \\
0 & \mbox{($\chi$ : odd, $\chi \neq \omega$)} \\
d_\chi \deg F(T) & \mbox{($\chi$ : even ).}
\end{array} \right. \]
If $S_\chi$ is empty, then 
$\zprank X_S (K_\infty)_\chi =\zprank X (K_\infty)_\chi$.
}

\par \bigskip

\noindent
\textbf{Proof. }
We may assume that $S$ is not empty.
Recall the following exact sequence: 
\[ E_\infty \rightarrow \bigoplus_{q \in S} R_q \rightarrow 
\Gal (M_S (K_\infty)/L(K_\infty)) \rightarrow 0 \]
which is stated in section 2.
By Proposition 6.2, 
we see that $(R_q)_\chi \otimes_{\mathbb{Z}_p} \mathbb{Q}_p$ is non-trivial 
if and only if $q \in S_\chi$.
Hence, if $S_\chi$ is empty, then 
$\Gal (M_S (K_\infty)/L(K_\infty))_\chi \otimes_{\mathbb{Z}_p} \mathbb{Q}_p$ is trivial,
and $\zprank X_S (K_\infty)_\chi =\zprank X (K_\infty)_\chi$.
In the following, we assume that $S_\chi$ is not empty.
By taking the $\chi$-quotient and tensoring with $\mathbb{Q}_p$, we have the exact sequence: 
\[ (E_\infty)_\chi \otimes_{\mathbb{Z}_p} \mathbb{Q}_p 
\overset{\eta_{\chi}}{\longrightarrow} 
\bigoplus_{q \in S_\chi} (R_q)_\chi \otimes_{\mathbb{Z}_p} \mathbb{Q}_p \rightarrow 
\Gal (M_S (K_\infty)/L(K_\infty))_\chi \otimes_{\mathbb{Z}_p} \mathbb{Q}_p \rightarrow 0. \]
It is sufficient to determine the cokernel of $\eta_{\chi}$.

Since $p$ is odd, we have a decomposition $E_\infty \cong E_\infty^+ \oplus E_\infty^-$, 
and we can see $E_\infty^- \cong \varprojlim \mu_{p^n}$.
We also note that $(E_\infty)_\chi \cong (E_\infty^+)_\chi \oplus (E_\infty^-)_\chi$
for a character $\chi$ of $G$.
It was already shown that 
if $\chi$ is an odd character, then $(E_\infty^+)_\chi$ is trivial,
and hence $(E_\infty)_\chi \cong (\varprojlim \mu_{p^n})_\chi$.

Assume that $\chi=\omega$.
Then $(E_\infty)_\omega \cong \varprojlim \mu_{p^n}$, and the natural mapping 
$\varprojlim \mu_{p^n} \rightarrow \bigoplus_{q \in S_\omega} R_q$ is injective.
We also note that $d_\omega=1$.
From these facts, we see that 
\[ \zprank X_S (K_\infty)_\omega
= \zprank X (K_\infty)_\omega + \sum_{q \in S_\omega} \zprank (R_q)_\omega -1 
= \zprank X (K_\infty)_\omega + \sum_{q \in S_\omega} p^{m_q} -1. \]

Assume that $\chi$ is odd and $\chi \neq \omega$.
Then $(E_\infty)_\chi \otimes_{\mathbb{Z}_p} \mathbb{Q}_p$ is trivial.
Hence we see 
\[ \zprank X_S (K_\infty)_\chi 
= \zprank X (K_\infty)_\chi + \sum_{q \in S_\chi} d_\chi p^{m_q}. \]

Let $\varepsilon$ be the trivial character.
Then we can see that 
\[ X_S (K_\infty)_{\varepsilon} \otimes_{\mathbb{Z}_p} \mathbb{Q}_p \cong 
X_S (\mathbb{Q}_\infty) \otimes_{\mathbb{Z}_p} \mathbb{Q}_p, \quad
X (K_\infty)_{\varepsilon} \otimes_{\mathbb{Z}_p} \mathbb{Q}_p \cong 
X (\mathbb{Q}_\infty) \otimes_{\mathbb{Z}_p} \mathbb{Q}_p. \]
Hence $\zprank X (K_\infty)_{\varepsilon} =0$.
We also note that 
\[ S_{\varepsilon} =\{ q \in S \; | \; \omega (\sigma_{0,q}) =1 \}
= \{ q \in S \; | \; q \equiv 1 \pmod{p} \}. \]
By the results for $\mathbb{Q}_\infty$ (see \cite{IMO}),
we see that $X_S (\mathbb{Q}_\infty) = X_{S_{\varepsilon}} (\mathbb{Q}_\infty)$, and
\[ \zprank X_S (K_\infty)_{\varepsilon} 
= \sum_{q \in S_{\varepsilon}} p^{m_q} - 
\mathrm{max} \{ p^{m_q} \; | \; q \in S_{\varepsilon} \}. \]
Since $d_{\varepsilon}=1$ and 
$f_{q, \varepsilon} (T) = (1+T)^{p^{m_q}}- \kappa_0^{p^{m_q}}$,
we see that 
\[ \mathrm{max} \{ p^{m_q} \; | \; q \in S_{\varepsilon} \}
= d_{\varepsilon} \mathrm{deg} F(T). \]
From these facts, the formula
\[ \zprank X_S (K_\infty)_{\varepsilon} = X (K_\infty)_{\varepsilon} 
+ \sum_{q \in S_{\varepsilon}} d_{\varepsilon} p^{m_q} - 
d_{\varepsilon} \mathrm{deg} F(T) \]
is certainly satisfied.

Finally, assume that $\chi$ is non-trivial and even.
By Proposition 6.3, we see that
\[ (E_\infty)_\chi /f_{q,\chi} (T) \otimes_{\mathbb{Z}_p} \mathbb{Q}_p
\rightarrow (R_q)_\chi \otimes_{\mathbb{Z}_p} \mathbb{Q}_p \]
is an isomorphism.
This isomorphism implies that 
\[ (E_\infty)_\chi /F(T) \otimes_{\mathbb{Z}_p} \mathbb{Q}_p
\rightarrow \bigoplus_{q \in S_\chi} (R_q)_\chi \otimes_{\mathbb{Z}_p} \mathbb{Q}_p \]
is injective.
By using the same argument stated in the proof of Proposition 6.3,
we obtain that 
$\dim_{\mathbb{Q}_p}
((E_\infty)_\chi /F(T)) \otimes_{\mathbb{Z}_p} \mathbb{Q}_p
= d_\chi \deg F(T)$.
Hence we see that 
\[ \zprank X_S (K_\infty)_\chi 
= \zprank X (K_\infty)_\chi + \sum_{q \in S_\chi} d_\chi p^{m_q}
- d_\chi \mathrm{deg} F(T). \]

We have shown the formula for all cases.
\hfill $\Box$

\par \bigskip

\noindent
\textbf{Remark. }
Assume that $p$ does not divide $|G|$.
Then $f_{q,\chi} (T)= (1+T)^{p^{m_q}}-\kappa_0^{p^{m_q}}$, and hence 
$\mathrm{deg} F(T) = \mathrm{max} \{ p^{m_q} \; | \; q \in S_\chi\}$.
Moreover, we can see that $p^{m_q}$ is equal to 
the number of primes of $\mathbb{Q}_\infty$ lying above $q$.

\par \bigskip

\noindent
\textbf{Example 6.5. }
We put $K=\mathbb{Q} (\mu_p)$, the $p$th cyclotomic field
(recall that $p$ is an odd prime).
In this case, every character of $G=\Gal (K/\mathbb{Q})$ is written by the form 
$\omega^i$ with $0 \leq i \leq p-2$.
Note also that $q \in S$ is unramified in $K$.
We may identify $G$ with $(\mathbb{Z}/p \mathbb{Z})^{\times}$, 
and $\sigma_{0,q}$ with $q \pmod{p}$. 
Then we can see 
\[ S_{\omega^i} = \{ q \in S \; | \; \omega^{1-i} (\sigma_{0,q})=1 \}
= \{ q \in S \; | \; i \equiv 1 \pmod{f_q} \}, \]
where $f_q$ is the order of $q$ in $(\mathbb{Z}/p \mathbb{Z})^{\times}$.
Assume that $S_{\omega^i}$ is not empty.
We put 
\[ P^{(i)} = \left\{ \begin{array}{ll}
1 & \mbox{($i=1$)} \\
0 & \mbox{($i$ : odd, $i \neq 1$)} \\
\max \{ p^{m_q} \; | \; q \in S_{\omega^i} \} & \mbox{($i$ : even)}. 
\end{array} \right. \]
By Theorem 6.4, we see
\[ \zprank X_S (K_\infty)_{\omega^i} = \lambda_{\omega^i} +
\sum_{q \in S_{\omega^i}} p^{m_q} -P^{(i)}, \]
where $\lambda_{\omega^i} = \zprank X (K_\infty)_{\omega^i}$ is the 
$\omega^i$-part of the (unramified)
Iwasawa $\lambda$-invariant of $K_\infty/K$.

\par \bigskip

\noindent
\textbf{Example 6.6. }
Let $k$ be a real quadratic field with conductor $d$
(the case that $p$ divides $d$ is allowed).
We put $K=k (\mu_p)$. 
Let $\chi$ be the quadratic character of $G=\Gal(K/\mathbb{Q})$ corresponding to $k$.
We may regard $\chi$ (resp. $\omega$) as a Dirichlet character
modulo $d$ (resp. modulo $p$).
In this case, we see
\[ S_\chi = \{ q \in S \; | \; \chi (q) \neq 0, \;
\chi (q) = \omega (q) \}. \]
Hence $S_\chi$ consists of the primes in $S$ which satisfy:
\begin{itemize}
\item $q \equiv 1 \pmod{p}$ and $q$ splits in $k$, or
\item $q \equiv -1 \pmod{p}$ and $q$ is inert in $k$.
\end{itemize}
Assume that $S_\chi \neq \emptyset$.
We put $P=\max \{ p^{m_q} \; | \; q \in S_\chi \}$, 
then we obtain the formula 
\[ \zprank X_S (K_\infty)_\chi = \zprank X (K_\infty)_\chi +
\sum_{q \in S_\chi} p^{m_q} -P. \]
Note that $\zprank X (K_\infty)_\chi$ is equal to the (unramified) Iwasawa $\lambda$-invariant 
of $k_\infty/k$. 
(If Greenberg's conjecture is true for $k$ and $p$, then 
$\zprank X (K_\infty)_\chi = 0$.)
Since 
\[ \zprank X_S (k_\infty) = 
\zprank X_S (K_\infty)_\chi + \zprank X_{S_{\varepsilon}} (\mathbb{Q}_\infty) \]
(where $\varepsilon$ is the trivial character),
we can obtain a formula of the $\mathbb{Z}_p$-rank of $X_S (k_\infty)$
(including $\zprank X (K_\infty)_\chi$).

\par \bigskip

\noindent
\textbf{Example 6.7. }
Let $F/\mathbb{Q}$ be a cyclic extension of degree $p$. 
Assume that $p$ is unramified in $F$.
We put $K=F (\mu_p)$, and
fix a character $\chi$ of $G=\Gal(K/\mathbb{Q})$ satisfying $K^{\mathrm{ker} (\chi)} =F$.
Let $\sigma$ be a fixed generator of 
$\Gal(F_\infty/\mathbb{Q}_\infty) \cong \Gal(F/\mathbb{Q})$, 
and $F^{(i)}$ the fixed field of $F_\infty$ by 
$\overline{\langle \gamma \sigma^i \rangle}$ for $0 \leq i \leq p-1$
(hence $F^{(0)}=F$).
For simplicity, we assume that every prime of $S$ is not decomposed in $\mathbb{Q}_\infty$.
Hence if $q \in S$ is unramified in $F$, then 
the splitting field of $F_\infty/\mathbb{Q}$ for $q$ must be one of 
$F^{(0)}, \ldots$, or $F^{(p-1)}$.
By the definition of $\chi$, we see that $\chi^{-1} (\sigma)$ is defined, 
and we put $\chi^{-1} (\sigma) =\zeta$ 
(note that $\zeta$ is a primitive $p$th root of unity).
In this case, we obtain that
\[ S_\chi = \{ q \in S \; | \; q \equiv 1 \pmod{p}, \;
\mbox{$q$ is not ramified in $F$} \}. \]
Under the assumption for $S$, we see $m_q=0$ for all $q \in S_\chi$.
When $q \in S_\chi$ splits in $F^{(i)}$, 
we see that $\chi^{-1} (\sigma_{p,q}) = \chi^{-1} (\sigma^i) = \zeta^i$,
and hence $f_{q,\chi} (T) = (1+T)- \zeta^i \kappa_0$.
We note that if $i \neq j$, then 
$(1+T)- \zeta^i \kappa_0$ and $(1+T)- \zeta^j \kappa_0$ are 
relatively prime.
We put 
\[ S_{\chi,i} = \{ q \in S_\chi \; | \; 
\mbox{$q$ splits in $F^{(i)}$} \}, \]
for $0 \leq i \leq p-1$.
Assume that $S_{\chi} \neq \emptyset$.
From the above facts, we see that 
$\deg F(T)$ is equal to the number of non-empty $S_{\chi, i}$'s.
That is, 
\[ \deg F(T)= |\Psi|, \; \mbox{where }
\Psi = \{ i \; |\; 0 \leq i \leq p-1, \; S_{\chi, i} \neq \emptyset\}. \]
Since $d_\chi=p-1$, we see 
\[ \zprank X_S (K_\infty)_\chi 
= \zprank X (K_\infty)_\chi + (p-1) \sum_{i \in \Psi}
\left( |S_{\chi,i}| -1 \right). \]

\par \bigskip

\noindent
\textbf{Acknowledgements. } 
The central part of the idea in the proof of our main theorem 
is due to the previous paper \cite{IMO}.
The author express his thanks to Professors Yasushi Mizusawa and Manabu Ozaki
for giving many suggestions during (and after) the collaborative work.
The author also express his thanks to Dr.~Takae Tsuji for giving comments.
The author was supported by Research Grant of Research Institute of Chiba Institute of 
Technology.

\begin{flushleft}
Tsuyoshi Itoh \\
Division of Mathematics, 
Education Center, 
Faculty of Social Systems Science, 
Chiba Institute of Technology, 
2-1-1 Shibazono, Narashino, Chiba 275-0023, Japan \\
e-mail : \texttt{tsuyoshi.itoh@it-chiba.ac.jp}
\end{flushleft}

\end{document}